\documentclass[twoside,11pt]{article}
\title{} \author{} \date{}
\usepackage{latexsym,amssymb,times}
\input amssym.def
\newtheorem{te}{Theorem}[section]
\newtheorem{prop}[te]{Proposition}

\newtheorem{cor}[te]{Corollary}
\newtheorem{fac}[te]{Fact}
\newtheorem{lem}[te]{Lemma}
\newtheorem{cla}[te]{Claim}

\newtheorem{ex}[te]{Example}

\def\dok{\noindent{\bf Proof. }}
\def\kdok{\hfill $\Box$ \par \vspace*{2mm} }
\def\a{\alpha}
\def\b{\beta}
\def\g{\gamma}
\def\f{\varphi}
\def\p{\psi}
\def\o{\omega}
\def\k{\kappa}

\def\r{\rho}
\def\s{\sigma}

\def\t{\tau}
\def\S{{\mathbb S}}
\def\P{{\mathbb P}}
\def\Q{{\mathbb Q}}

\def\N{{\mathbb N}}
\def\X{{\mathbb X}}
\def\Y{{\mathbb Y}}
\def\Z{{\mathbb Z}}

\def\CB{{\mathcal B}}
\def\ED{{\mathcal E}{\mathcal D}_{\mathrm{fin}}}

\def\I{{\mathcal I}}

\def\U{{\mathcal U}}

\def\down{\!\downarrow}
\def\up{\!\uparrow}
\def\la{\langle}
\def\ra{\rangle}

\def\dom{\mathop{\mathrm{dom}}\nolimits}

\def\Lim{\mathop{\mbox{Lim}}\nolimits}

\def\Emb{\mathop{\rm Emb}\nolimits}

\def\At{\mathop{\rm At}\nolimits}

\def\Fin{\mathop{\rm Fin}\nolimits}

\def\rp{\mathop{\rm rp}\nolimits}

\def\sm{\mathop{\rm sm}\nolimits}
\def\sq{\mathop{\rm sq}\nolimits}

\def\ar{\mathop{\rm ar}\nolimits}
\begin{document}
\thispagestyle{plain}
\begin{center}
           {\Large \bf {From $\bf A_1$ to $\bf D_5$: 
           Towards a Forcing-Related Classification\\[2mm] of Relational Structures}}
\end{center}
\begin{center}
{\small\bf Milo\v s S.\ Kurili\'c}\\[1mm]
        {\small Department of Mathematics and Informatics, University of Novi Sad, \\
         Trg Dositeja Obradovi\'ca 4, 21000 Novi Sad, Serbia.\\[-1mm]
                                     e-mail: milos@dmi.uns.ac.rs}
\end{center}

\begin{abstract}
\noindent
We investigate the partial orderings of the form $\la \P (\X ), \subset\ra$, where $\X $ is a relational structure and
$\P (\X )$ the set of the domains of its isomorphic substructures.
A rough classification of countable binary structures
corresponding to the forcing-related properties of the posets of their copies is obtained.

{\sl 2000 Mathematics Subject Classification}:
03C15,  
03E40,  
06A10.  

{\sl Keywords}: relational structure, isomorphic substructure, poset,  forcing.
\end{abstract}

\section{Introduction}

The relational structure $\X= \la \o , <\ra$, where $<$ is the natural order on the set $\o$ of natural numbers
is a structure having the following extremal property: each $\o$-sized subset $A$ of $\o$ determines a substructure
isomorphic to the whole structure. If instead of $\la \o , <\ra$ we take the integer line $\Z =\la Z, <\ra$, then
we lose the maximality of the set of isomorphic substructures (the set of positive integers is not a copy of $\Z$). Finally, the minimality of the set of copies is reached by the linear graph $G_\Z =\la Z, \r\ra$, where $\r = \{ \la m, n \ra : |m-n|=1 \}$, since each proper subset
$A$ of $Z$ determines a disconnected graph and, hence, fails to be a copy of the whole graph.

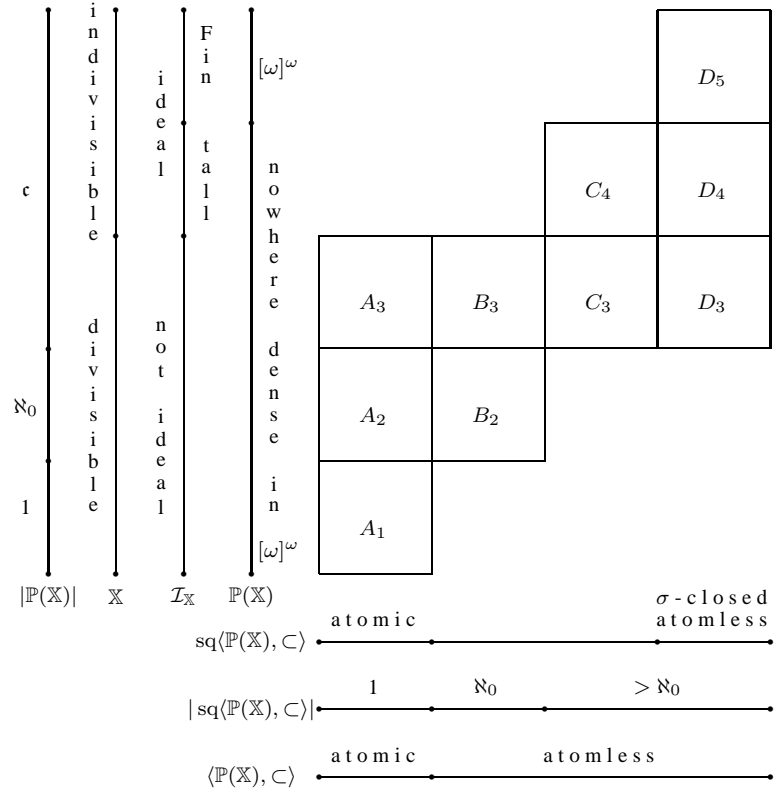
\begin{figure}[htb]\label{FIG2}
\begin{center}
\unitlength 0.6mm 
\linethickness{0.4pt}
\ifx\plotpoint\undefined\newsavebox{\plotpoint}\fi 


\begin{picture}(180,190)(0,0)


\put(70,10){\line(1,0){100}}
\put(70,25){\line(1,0){100}}
\put(70,40){\line(1,0){100}}
\put(10,55){\line(0,1){125}}
\put(25,55){\line(0,1){125}}
\put(40,55){\line(0,1){125}}
\put(55,55){\line(0,1){125}}
\put(70,55){\line(0,1){75}}
\put(95,55){\line(0,1){75}}
\put(120,80){\line(0,1){50}}
\put(120,130){\line(0,1){25}}
\put(145,105){\line(0,1){75}}
\put(170,105){\line(0,1){75}}
\put(70,55){\line(1,0){25}}
\put(70,80){\line(1,0){50}}
\put(70,105){\line(1,0){100}}
\put(70,130){\line(1,0){100}}
\put(120,155){\line(1,0){50}}
\put(145,180){\line(1,0){25}}


\put(70,10){\circle*{1}}
\put(95,10){\circle*{1}}
\put(170,10){\circle*{1}}
\put(70,25){\circle*{1}}
\put(95,25){\circle*{1}}
\put(120,25){\circle*{1}}
\put(170,25){\circle*{1}}
\put(70,40){\circle*{1}}
\put(145,40){\circle*{1}}
\put(170,40){\circle*{1}}
\put(10,55){\circle*{1}}
\put(10,80){\circle*{1}}
\put(10,105){\circle*{1}}
\put(10,180){\circle*{1}}
\put(25,55){\circle*{1}}
\put(25,130){\circle*{1}}
\put(25,180){\circle*{1}}
\put(40,55){\circle*{1}}
\put(40,130){\circle*{1}}
\put(40,155){\circle*{1}}
\put(40,180){\circle*{1}}
\put(55,55){\circle*{1}}
\put(55,155){\circle*{1}}
\put(55,180){\circle*{1}}
\put(95,40){\circle*{1}}

\scriptsize
\put(82,15){\makebox(0,0)[cc]{a t o m i c}}%
\put(132,15){\makebox(0,0)[cc]{a t o m l e s s}}%
\put(82,45){\makebox(0,0)[cc]{a t o m i c}}%

\put(82,30){\makebox(0,0)[cc]{1}}%
\put(107,30){\makebox(0,0)[cc]{$\aleph _0$}}%
\put(145,30){\makebox(0,0)[cc]{$>\aleph _0$}}%

\put(157,50){\makebox(0,0)[cc]{$\sigma$ - c l o s e d}}%
\put(157,45){\makebox(0,0)[cc]{a t o m l e s s}}%

\put(5,70){\makebox(0,0)[cc]{1}}%
\put(5,92){\makebox(0,0)[cc]{$\aleph _0$}}%
\put(5,140){\makebox(0,0)[cc]{${\mathfrak c}$}}%
\put(20,180){\makebox(0,0)[cc]{i}}%
\put(20,175){\makebox(0,0)[cc]{n}}%
\put(20,170){\makebox(0,0)[cc]{d}}%
\put(20,165){\makebox(0,0)[cc]{i}}%
\put(20,160){\makebox(0,0)[cc]{v}}%
\put(20,155){\makebox(0,0)[cc]{i}}%
\put(20,150){\makebox(0,0)[cc]{s}}%
\put(20,145){\makebox(0,0)[cc]{i}}%
\put(20,140){\makebox(0,0)[cc]{b}}%
\put(20,135){\makebox(0,0)[cc]{l}}%
\put(20,130){\makebox(0,0)[cc]{e}}%
\put(20,110){\makebox(0,0)[cc]{d}}%
\put(20,105){\makebox(0,0)[cc]{i}}%
\put(20,100){\makebox(0,0)[cc]{v}}%
\put(20,95){\makebox(0,0)[cc]{i}}%
\put(20,90){\makebox(0,0)[cc]{s}}%
\put(20,85){\makebox(0,0)[cc]{i}}%
\put(20,80){\makebox(0,0)[cc]{b}}%
\put(20,75){\makebox(0,0)[cc]{l}}%
\put(20,70){\makebox(0,0)[cc]{e}}%
\put(35,165){\makebox(0,0)[cc]{i}}%
\put(35,160){\makebox(0,0)[cc]{d}}%
\put(35,155){\makebox(0,0)[cc]{e}}%
\put(35,150){\makebox(0,0)[cc]{a}}%
\put(35,145){\makebox(0,0)[cc]{l}}%
\put(35,110){\makebox(0,0)[cc]{n}}%
\put(35,105){\makebox(0,0)[cc]{o}}%
\put(35,100){\makebox(0,0)[cc]{t}}%
\put(35,90){\makebox(0,0)[cc]{i}}%
\put(35,85){\makebox(0,0)[cc]{d}}%
\put(35,80){\makebox(0,0)[cc]{e}}%
\put(35,75){\makebox(0,0)[cc]{a}}%
\put(35,70){\makebox(0,0)[cc]{l}}%
\put(45,175){\makebox(0,0)[cc]{F}}%
\put(45,170){\makebox(0,0)[cc]{i}}%
\put(45,165){\makebox(0,0)[cc]{n}}%
\put(45,150){\makebox(0,0)[cc]{t}}%
\put(45,145){\makebox(0,0)[cc]{a}}%
\put(45,140){\makebox(0,0)[cc]{l}}%
\put(45,135){\makebox(0,0)[cc]{l}}%
\put(61,167){\makebox(0,0)[cc]{$[\omega ]^{\omega }$}}%
\put(60,145){\makebox(0,0)[cc]{n}}%
\put(60,140){\makebox(0,0)[cc]{o}}%
\put(60,135){\makebox(0,0)[cc]{w}}%
\put(60,130){\makebox(0,0)[cc]{h}}%
\put(60,125){\makebox(0,0)[cc]{e}}%
\put(60,120){\makebox(0,0)[cc]{r}}%
\put(60,115){\makebox(0,0)[cc]{e}}%
\put(60,105){\makebox(0,0)[cc]{d}}%
\put(60,100){\makebox(0,0)[cc]{e}}%
\put(60,95){\makebox(0,0)[cc]{n}}%
\put(60,90){\makebox(0,0)[cc]{s}}%
\put(60,85){\makebox(0,0)[cc]{e}}%
\put(60,75){\makebox(0,0)[cc]{i}}%
\put(60,70){\makebox(0,0)[cc]{n}}%
\put(61,60){\makebox(0,0)[cc]{$[\omega ]^{\omega }$}}%



\put(10,50){\makebox(0,0)[cc]{$|\P (\X )|$}}%
\put(25,50){\makebox(0,0)[cc]{$\X$}}%
\put(40,50){\makebox(0,0)[cc]{${\mathcal I}_\X$}}%
\put(55,50){\makebox(0,0)[cc]{$\P(\X)$}}%
\put(55,40){\makebox(0,0)[cc]{$\sq \langle \P (\X ) , \subset \rangle$}}%
\put(55,25){\makebox(0,0)[cc]{$|\sq \langle \P (\X ) , \subset \rangle|$}}%
\put(55,10){\makebox(0,0)[cc]{$\langle \P (\X ) , \subset \rangle$}}%


\put(82,65){\makebox(0,0)[cc]{$A_1$}}%
\put(82,90){\makebox(0,0)[cc]{$A_2$}}%
\put(82,115){\makebox(0,0)[cc]{$A_3$}}%
\put(107,90){\makebox(0,0)[cc]{$B_2$}}%
\put(107,115){\makebox(0,0)[cc]{$B_3$}}%
\put(132,115){\makebox(0,0)[cc]{$C_3$}}%
\put(132,140){\makebox(0,0)[cc]{$C_4$}}%
\put(157,115){\makebox(0,0)[cc]{$D_3$}}%
\put(157,140){\makebox(0,0)[cc]{$D_4$}}%
\put(157,165){\makebox(0,0)[cc]{$D_5$}}%
\end{picture}

\end{center}

\vspace{-6mm}

\caption{Binary relations on countable sets}
\end{figure}

We investigate the posets of the form $\la \P (\X ), \subset\ra$, where $\X $ is a relational structure and
$\P (\X )$ the set of the domains of its isomorphic substructures.
Although some our statements are general,
the main result of the paper is the diagram on Figure \ref{FIG2},
describing an interplay between
the properties of a countable binary structure $\X$
and the properties of the corresponding poset $\la \P (\X ), \subset \ra$.
So we obtain a rough classification of countable binary structures
concerning the forcing-related properties of the posets of their copies:
for the structures from column A (resp.\ B; D) the corresponding posets are forcing equivalent to the trivial poset
(resp.\ the Cohen forcing, $\la {}^{<\o }2, \supset\ra$;
a $\s$-closed atomless poset) and the wild animals are in cages $C_3$ and $C_4$, where the posets of copies are forcing equivalent to the quotients
of the form $P(\o )/\I$, for some co-analytic tall ideal $\I$.

Clearly, such classification
depends on the model of
set theory in which we work. For example, under the CH all the structures from column $D$ are in the same class
(having the posets of copies forcing equivalent to $(P(\o)/\Fin)^+$), but this is not true
in the Mathias model. Also the classification is very rough. Namely, it is easy to see that equimorphic structures have forcing equivalent posets of copies \cite{Kurstr} and, hence,
all countable non-scattered linear orders are equivalent in this sense. Moreover,
the class of structures satisfying $\P (\X )=\{ X \}$ contains continuum many non-equimorphic structures \cite{Kemin}.

\noindent
A few words on notation.
Let $L=\{ R_i :i\in I \}$ be a relational language, where $\ar (R_i)=n_i$, $i\in I$.
An $L$-structure $\X = \la X, \{ \r _i :i\in I \} \ra$ is called {\it countable} iff $|X|=\o$; {\it binary} iff $L=\{ R \}$ and $\ar (R)=2$.
If  $A\subset X$, then $\la A, \{ (\r _i)_A :i\in I \} \ra$ is a {\it substructure}
of $\X$, where $(\r_i)_A = \r _i \cap A^{n_i}$, $i\in I$.
If $\Y =\la Y, \{ \s _i : i\in I\}  \ra$ is an $L$-structure too,
a mapping $f:X \rightarrow Y$ is an {\it embedding} (we write
$\X \hookrightarrow _f \Y$) iff it is an injection and
$$\forall i\in I \;\;\forall \la x_1, \dots   x_{n _i}\ra \in X^{n_i} \;\;
(\la x_1, \dots   ,x_{n _i}\ra \in \r _i
\Leftrightarrow
\la f(x_1), \dots  ,f(x_{n _i})\ra \in \s _i).
$$
If $\X$ embeds in $\Y$ we write $\X \hookrightarrow \Y$. Let
$\Emb (\X , \Y )  =  \{ f: \X \hookrightarrow _f \Y \} $ and
$\Emb (\X )  = \{ f: \X \hookrightarrow _f \X \}$.
If, in addition, $f$ is a surjection, it is an {\it isomorphism} (we write
$\X \cong _f \Y$) and the structures $\X$ and $\Y$ are {\it isomorphic},
in notation $\X \cong \Y$.
So we investigate the posets of the form $\la \P (\X ), \subset \ra$, where $\X=\la X, \{ \r _i :i\in I \} \ra$ is a relational structure and
$$
\P (\X )  =  \{ A\subset X : \la A, \{ (\r _i)_A :i\in I \} \ra \cong \X\}
          =  \{ f[X] : f \in \Emb (\X )\} .
$$
More generally, if  $\Y =\la Y, \{ \s _i : i\in I\}  \ra$ is a structure of the same language, let
$\P (\X , \Y  )  =  \{ B\subset Y : \la B, \{ (\s _i)_B :i\in I \} \ra \cong \X \}
                =  \{ f[X] : f \in \Emb (\X , \Y )\} $.
\section{Homogeneity and atoms}

If $\P =\la P , \leq \ra$ is a partial order, $p,q\in \P $ are {\it compatible} iff there is $r\leq p,q$.
Otherwise $p$ and $q$ are {\it incompatible} and we write $p\perp q$.
$p\in P$ is an {\it atom}, in notation  $p\in \At (\P )$, iff each  $q,r\leq p$ are compatible.
$\P $ is called: {\it atomless} iff $\At (\P )=\emptyset$; {\it atomic} iff $\At (\P )$ is dense in $\P$;
{\it homogeneous} iff  it has the largest element and $\P \cong p\down =(-\infty ,p]_{\P }$, for each
$p\in P$. Clearly we have
\begin{fac}\rm\label{T4018}
A homogeneous poset $\P =\la P ,\leq \ra$ is either atomless or downwards directed and $\At (\P)=P$ in the second case.
\end{fac}
A family $\CB $ is an {\it uniform filter base} on a set $X$ iff
(UFB1) $\emptyset \neq \CB \subset [X]^{|X|}$;
(UFB2) For each $A,B\in \CB$ there is $C\in \CB$ such that $C\subset A\cap B$.
\begin{te}\rm\label{T4013}
Let $\X = \la X, \{ \r _i :i\in I \} \ra$ be a relational structure. Then

(a) $\la \P (\X ), \subset \ra$ is a homogeneous poset;

(b) $\la \P (\X ), \subset \ra$ is either atomless or atomic;

(c) $\la \P (\X ), \subset \ra$ is atomless iff it contains two incompatible elements;

(d) If $\la \P (\X ), \subset \ra$ is atomic, then $\At (\P (\X ) )=\P (\X )$ and, moreover, $\P (\X )$ is
an uniform filter base on $X$.  Also $\bigcap \P (\X )\in \P (\X )$ iff  $\P (\X )=\{ X\}$.
\end{te}
\dok
(a) Clearly, $1_{\P (\X )}=X$. Let $C\in \P (\X )$ and $f\in \Emb (\X )$, where $C=f[X]$.
We show that $\la \P (\X ), \subset \ra\cong _F \la  (-\infty , C]_{\P (\X )} , \subset \ra$,
where the function $F$ is defined by $F(A)=f[A]$, for each $A\in \P (\X )$. For $A\in \P (\X )$
we have $F(A)\subset C$ and there is $g\in \Emb (\X )$ such that $A=g[X]$. Clearly
$f\circ g \in \Emb (\X )$ and, hence, $F(A)=f[g[X]]\in \P(\X)$. Thus $F: \P (\X )\rightarrow  (-\infty , C]_{\P (\X )} $.

Since $f$ is an injection, $f[A]=f[B]$ implies $A=B$, so $F$ is an injection.

Let $\P (\X )\ni B \subset C$. Since $B\subset f[X]$ we have $B=f[f^{-1}[B]]$
and, clearly,  $\la f^{-1}[B] ,\{ (\r _i)_{f^{-1}[B]} :i\in I \} \ra \cong _{f|f^{-1}[B]} \la B, \{ (\r _i)_B :i\in I \} \ra\cong \X$.
Thus $f^{-1}[B]\in \P (\X )$  and $B=F(f^{-1}[B])$, so $F$ is a surjection.

Since $f$ is an injection, for $A,B\in \P (\X )$ we have $A\subset B \Leftrightarrow f[A]\subset f[B]$. Thus $F$ is an
order isomorphism.

(b) Follows from (a) and Fact \ref{T4018}.

(c) If $\P (\X)$ contains two incompatible elements, then it is not downwards directed and, by Fact \ref{T4018},
must be atomless.

(d) Let $\la \P (\X ), \subset \ra$ be atomic. By Fact \ref{T4018}, $\At (\P (\X ) )=\P (\X )$ and $\P (\X )$ satisfies
(UFB2). Since $X\in \P (\X )\subset [X]^{|X|}$, (UFB1) holds as well.
Suppose that $A= \bigcap \P (\X )\in \P (\X )$ and $\P (\X )\neq \{ X\}$. Then $A \varsubsetneq X$ and, since
$\P (\X )\cong A\down$, there is $B\in \P (\X )$ such that $B\varsubsetneq A$. A contradiction.
\hfill $\Box$
\section{The complexity and size}
For each relational structure $\X$ we have $\{ X \} \subset \P (\X ) \subset [X]^{|X|}$ and $\P (\X )$ is of size 1 or infinite, because
if $f\in \Emb (\X)$ and $f[X]\neq X$, then $f^n [X]$, $n\in\N$, is a decreasing sequence of elements of $\P (\X)$.
Now we show that $|\P (\X )|\in \{ 1, \aleph _0 ,{\mathfrak c}\}$.

By $2^\o$ and $\o ^\o$ we denote the Cantor cube and the Baire space and
$p_k:2^\o \rightarrow 2$ and $\pi _k : \o ^\o \rightarrow \o$, $k\in \o$, will be the corresponding projections.
As usual, the mapping $\chi : P(\o )\rightarrow 2^\o$, where
$\chi (A)=\chi _A$, for each $A\subset \o$, identifies the subsets of $\o $ with their characteristic functions
and a set ${\mathcal S}\subset P(\o )$ is called closed (Borel, analytic ...) iff
$\chi [{\mathcal S} ] $ is a closed (Borel, analytic ...) set in the space $2^\o$.

For ${\mathcal S}\subset P(\o )$ let ${\mathcal S}\up =\{ A\subset \o : \exists S\in {\mathcal S} \; S\subset A \}$ and, for $A\subset 2^ \o  $,
let $A\up =\{ x\in 2^ \o : \exists a\in A \; a\leq x \}$, where $a\leq x$ means that
$a(n)\leq x(n)$, for all $n\in \o$. Instead of $\{ a\} \up$ we will write $a\up$.
\begin{te}  \rm  \label{T4052}
If $\X =\la X , \{ \r _i :i\in I \}\ra$ is a countable relational structure and ${\mathcal I}_\X =\{ I\subset X : \neg \exists A \in \P (\X )\; A\subset I \}$, then

(a) $\P (\X)$  is an analytic set;

(b) $\P (\X) \up$ is an analytic set;

(c) ${\mathcal I}_\X $ is a co-analytic set containing the ideal $\Fin _X$ of finite subsets of $X$;

(d) The sets $\P (\X)$ and $\P (\X)\! \up$ have the Baire property and size 1, $\aleph _0 $ or ${\mathfrak c}$.
\end{te}
\dok
Without loss of generality  we suppose $X=\o $. Let $\ar (\r _i)=n_i$, $ i\in I$.

(a) This statement is a folklore but, for completeness, we include its proof.

\vspace{2mm}
\noindent
{\it Claim 1.}
$\Emb (\X )$ is a closed set in the Baire space, $\o ^\o$.

\vspace{2mm}
\noindent
{\it Proof of Claim 1.}
We show that the set $\o ^\o \setminus\Emb (\X )$ is open. Let $f\in \o ^\o \setminus\Emb (\X )$.

If $f$ is not an injection and $m,n\in \o$, where $m\neq n$ and $f(m)=f(n)=k$, then $\pi ^{-1}_m [\{ k\}] \cap \pi ^{-1}_n [\{ k\}]$ is a
neighborhood of $f$ contained in  $\o ^\o \setminus\Emb (\X )$.

Otherwise there are $i\in I$ and $m_1, \dots , m_{n_i}\in \o$ such that
$\la m_1, \dots , m_{n_i} \ra \in \r _i \not\Leftrightarrow \la f(m_1), \dots ,f(m_{n_i}) \ra \in \r _i$.
Then $B=\bigcap _{j\leq n_1}\pi _{m_j}^{-1}[\{ f(m_j)\}]$ is a neighborhood of $f$ contained in  $\o ^\o \setminus\Emb (\X )$.

\vspace{2mm}
\noindent
{\it Claim 2.}
The mapping $F: \o ^\o \rightarrow 2^\o $ defined by $F(f)=\chi _{f[\o ]}$ is a Borel mapping.

\vspace{2mm}
\noindent
{\it Proof of Claim 2.}
By \cite{Kech}, p.\ 71, it is sufficient to show that $F^{-1}[p_n^{-1}[\{ j \}]]$ is a Borel set, for each $n\in \o$ and $j\in 2$.
Clearly, for $f\in \o ^\o$ we have   $f\in F^{-1}[p_n^{-1}[\{ j \}]] \Leftrightarrow \chi _{f[\o ]}(n)=j$. Thus $f\in F^{-1}[p_n^{-1}[\{ 1 \}]]$ iff
$n\in f[\o ]$ iff $f(k)=n$, that is $f\in \pi _k ^{-1}[\{ n\}]$, for some $k\in \o$. So $F^{-1}[p_n^{-1}[\{ 1 \}]] =\bigcup _{k\in \o}\pi _k ^{-1}[\{ n\}]$
is an open set and, similarly, $F^{-1}[p_n^{-1}[\{ 0 \}]] =\o ^\o \setminus \bigcup _{k\in \o}\pi _k ^{-1}[\{ n\}]$ is closed and, hence, Borel.

\vspace{2mm}
\noindent
{\it Claim 3.}
$\chi [\P (\X )]=F [\Emb (\X )]$.

\vspace{2mm}
\noindent
{\it Proof of Claim 3.}
Since $\chi$ is a bijection, for $A\subset \o$ we have: $\chi _A \in \chi [\P (\X )]$ iff $A\in \P (\X )$ iff $A=f[\o ]$, that is
$\chi _A= \chi _{f[\o ]}=F(f)$, for some $f\in \Emb (\X)$ iff $\chi _A \in F [\Emb (\X )]$.

\vspace{2mm}
\noindent
By Claims 1 and 2, $F [\Emb (\X )]$ is an analytic set (see e.g.\ \cite{Kech}, p.\ 86).
Thus, by Claim 3, the set $\chi [\P (\X )]$ is analytic.

(b) If we regard the set $\Emb (\X)$ as a subspace of the Baire space $\o ^\o$, then
$\{ \pi _k^{-1 }[\{ n \}] \cap \Emb (\X) : k,n\in \o \}$ is a subbase for the corresponding topology on $\Emb (\X)$ and we have

\vspace{2mm}
\noindent
{\it Claim 4.}
$B= \bigcup _{f\in \Emb (\X )} \{ f \} \times \chi _{f[\o ]} \up $ is a closed set in the product $\Emb (\X )\times 2^\o$.

\vspace{2mm}
\noindent
{\it Proof of Claim 4.}
Let $\la f,x \ra \in (\Emb (\X )\times 2^\o )\setminus B$. Then $x\not\in \chi _{f[\o ]} \up $ and, hence, there is $n_0\in \o $ such that
$x(n_0)<\chi _{f[\o ]} (n_0)$. Thus, first, $x(n_0)=0$, which implies
$x\in p_{n_0}^{-1}[\{ 0 \}]$
and, second, $\chi _{f[\o ]} (n_0)=1$, that is $n_0\in f[\o]$ so there is $k_0\in \o$ satisfying $f(k_0)=n_0$ and, hence,
$f\in \pi_{k_0}^{-1}[\{ n_0 \}]$.
Now we have $\la f,x \ra \in O =(\pi_{k_0}^{-1}[\{ n_0 \}] \cap \Emb (\X ) ) \times p_{n_0}^{-1}[\{ 0 \}]$
and we show that $O\cap B=\emptyset$. Suppose that $\la g,y \ra \in O\cap B$. Then, since  $\la g,y \ra \in O$, we have
$g(k_0)=n_0$  and $ y(n_0)=0$; since $\la g,y \ra \in B$ we have
$y\geq \chi _{g[\o ]}$, which implies
$\forall n\in g[\o ] \;\; y(n)=1$. So $ y(n_0)=0$ implies $n_0\not\in g[\o ]$, which is not true because $g(k_0)=n_0$.
Thus $O$ is a neighborhood of $\la f,x \ra$ contained in $(\Emb (\X )\times 2^\o )\setminus B$ and this set is open.

\vspace{2mm}
\noindent
{\it Claim 5.}
$\chi [\P (\X )\up ] = \pi _{2^\o }[B]$, where $\pi _{2^\o }: \Emb (\X )\times 2^\o\rightarrow 2^\o$ is the projection.

\vspace{2mm}
\noindent
{\it Proof of Claim 5.}
If $x\in\chi [\P (\X )\up ] $, then there are $C\in \P (\X )$ and $A$ such that $C\subset A \subset \o$ and
$x=\chi _A$. Let $f\in \Emb (\X)$, where $C=f[\o ]$. Then $f[\o ]\subset A$ implies $x\geq \chi _{f[\o ]}$
and, hence, $\la f,x \ra \in B$ and $x =  \pi _{2^\o }(\la f,x \ra ) \in \pi _{2^\o }[B]$.

If  $x \in \pi _{2^\o }[B]$, then there is $f\in \Emb (\X)$  such that $x\geq \chi _{f[\o ]}$ and for $A= x^{-1}[\{ 1 \}]$
we have $x=\chi _A \geq \chi _{f[\o ]}$, which implies $\P (\X )\ni f[\o ]\subset A$, that is $A\in \P (\X )\up$ and, hence,
$x=\chi (A)\in \chi [\P (\X )\up ]$.

\vspace{2mm}
\noindent
By Claim 1, $\Emb (\X )$ is a Polish space so $\Emb (\X )\times 2^\o$ is a Polish space too.
Since the projection $\pi _{2^\o }$ is continuous, it is a Borel mapping and, by Claim 4,
$\pi _{2^\o }[B]$ is an analytic set (see \cite{Kech}, p.\ 86). By Claim 5 the set  $\chi [\P (\X )\up ]$ is analytic as well.

(c) follows from (b) and the equality $\I _\X = P(X)\setminus \P (\X )\up$.

(d) follows from (a), (b) and known facts about analytic sets (see \cite{Kech}).
\hfill $\Box$
\section{The separative quotient}
A partial order $\P =\la P , \leq \ra$ is called
{\it separative} iff for each $p,q\in P$ satisfying $p\not\leq q$ there is $r\in P$ such that $r\leq p$ and $r\perp q$.
The {\it separative modification} of $\P$
is the separative pre-order $\sm (\P )=\la P , \leq ^*\ra$, where
$p\leq ^* q$ iff $\forall r\leq p \; \exists s \leq r \; s\leq q $.
The {\it separative quotient} of $\P$
is the separative partial order $\sq (\P )= \langle P /\!\! =^* , \trianglelefteq \rangle$, where
$p = ^* q \Leftrightarrow p \leq ^* q \land q \leq ^* p\;$  and $\;[p] \trianglelefteq [q] \Leftrightarrow p \leq ^* q $.

If $\k$ is a regular cardinal, a pre-order $\P =\langle P , \leq \rangle$ is $\k${\it -closed} iff for each
$\gamma <\k$ each sequence $\la p_\a :\a <\gamma\ra$ in $P$, such that $\a <\beta \Rightarrow p_{\beta}\leq p_\a$,
has a lower bound. $\o _1$-closed pre-orders are called $\s${\it -closed}
and the following facts are well known.
\begin{fac} \rm \label{T4055}
Let $\P $ be a partial order. Then

(a) $\P$, $\sm (\P)$ and $\sq (\P)$
are forcing equivalent forcing notions;

(b) $\P$ is atomless iff
$\sm (\P )$ is atomless iff $\sq (\P )$ is atomless.
\end{fac}
\begin{fac} \rm \label{T4056}
If $\kappa ^{<\kappa }=\kappa$, then all atomless separative $\kappa$-closed pre-orders  of size $\kappa$, are
forcing equivalent (for example to the tree $\la {}^{<\kappa}\kappa , \supset\ra$).
\end{fac}
\begin{te}\rm \label{T4057}
Let $\X =\la X , \{ \r_i :i\in I \}\ra$  be a relational structure. Then

(a) $\sm \la \P (\X ), \subset \ra = \la \P (\X ), \leq ^* \ra $, where for $A, B \in \P (\X )$
\begin{equation}\label{EQ4002}
A\leq ^* B \Leftrightarrow \forall C \in \P (\X )\;\;
(C\subset A \Rightarrow \exists D\in \P (\X ) \;\; D\subset C\cap B);
\end{equation}

(b) $|\sq \la \P (\X ), \subset \ra |=1$ iff $\la \P (\X ), \subset \ra$ is atomic;

(c) $|\sq \la \P (\X ), \subset \ra |\geq \aleph _0$ iff $\la \P (\X ), \subset \ra$ is atomless;

(d) If $|\sq \la \P (\X ), \subset \ra | = \aleph _0$, then $\la \P (\X ), \subset \ra$ is forcing equivalent to the reversed binary tree
$\la {}^{<\o }2, \supset \ra$ (a forcing notion adding one Cohen real);

(e) If CH holds and $\sq \la \P (\X ), \subset \ra $ is $\s$-closed, atomless and of size $\mathfrak c$, then $\la \P (\X ), \subset \ra$ is forcing equivalent
to $(P(\o)/\Fin)^+$.
\end{te}
\dok
(a) This follows directly from the definition of the separative modification.

(b) If $|\sq \la \P (\X ), \subset \ra |=1$, then for each $A,B \in \P (\X )$ we have $A\leq ^* B$ so, by (\ref{EQ4002}), there is
$D\in \P (\X )$ such that $D\subset A\cap B$. Thus $\la \P (\X ), \subset \ra$ is downwards directed and, hence, atomic.

If $\la \P (\X ), \subset \ra$ is atomic and $A,B \in \P (\X )$, then, by Theorem \ref{T4013}(d), for each $C\in \P (\X)$ satisfying
$C\subset A$ there is $D\in \P (\X) $ such that $D\subset C\cap B$. Thus, by (\ref{EQ4002}), $A\leq ^* B$, for each $A, B\in \P (\X)$. Hence  $A= ^* B$, for each $A, B\in \P (\X)$,
and, consequently, $|\sq \la \P (\X ), \subset \ra |=1$.

(c) The implication ``$\Rightarrow$" follows from (b) and Theorem \ref{T4013}(b).
If the poset $\la \P (\X ), \subset \ra$ is atomless, then it contains an infinite antichain $\{ A_n : n\in \o\}$.
By (a), $A\leq ^* B$ implies that $A$ and $B$ are compatible,
thus $A_m \neq ^* A_n$, for $m\neq n$, which implies that the set $\sq \la \P (\X ), \subset \ra$ is infinite.

(d) If $|\sq \la \P (\X ), \subset \ra | = \aleph _0$, then, by (c), the partial order $\la \P (\X ), \subset \ra$ is atomless and, by Fact
\ref{T4055}(b), $\sq \la \P (\X ), \subset \ra$ is atomless as well. By Facts \ref{T4055}(a) and \ref{T4056} (for $\kappa =\o$), $\la \P (\X ), \subset \ra$ is forcing equivalent to the forcing $\la {}^{<\o }\o, \supset \ra$ or to $\la {}^{<\o }2, \supset \ra$.

(e) follows from Facts \ref{T4055}(a) and \ref{T4056} (for $\kappa =\o _1$).
\hfill $\Box$
\begin{ex}\rm \label{EX4014}
$\la \P (\X ),\subset \ra$ is a separative poset isomorphic to $\la {}^{<\o }2 , \supset \ra$.
Let $G_{{}^{<\o }2}$ be the digraph $\la {}^{<\o }2 , \r\ra$, where
$\r =\{ \la \f , \f ^\smallfrown i \ra : \f \in {}^{<\o }2 \land i\in 2 \}$.
For $\f \in {}^{<\o }2$ let $A_\f =\{\p \in {}^{<\o }2 : \f \subset \p \}$ and let us prove that
\begin{equation}\label{EQ4071}
 \P (G_{{}^{<\o }2} ) =\{ A_\f : \f \in {}^{<\o }2 \} .
\end{equation}
The inclusion ``$\supset$" is evident.
Conversely, if $A\in \P (G_{{}^{<\o }2} )$ and $f : G_{{}^{<\o }2}  \hookrightarrow G_{{}^{<\o }2}$, where $A=f[{}^{<\o }2]$, we show that
$A=A_{f(\emptyset )}$.

First, if $f (\f )\in A$ and $\dom (\f )=n$, then, since $\la \f \upharpoonright k , \f \upharpoonright (k+1) \ra\in \r$, for $k<n-1$,
we have  $\la f(\f \upharpoonright k ), f(\f \upharpoonright (k+1)) \ra\in \r$, for $k<n$. But this is an oriented path from
$f(\f \upharpoonright 0) =f(\emptyset )$ to  $f(\f \upharpoonright n) =f(\f )$, which implies
$f(\emptyset )\subset f(\f )$, that is $f(\f )\in A_{f(\emptyset )}$.
Second, by induction we show that $f(\emptyset )^\smallfrown \eta \in A$, for all $\eta \in {}^{<\o }2$. Let  $f(\emptyset )^\smallfrown \eta \in A$.
Then $f(\emptyset )^\smallfrown \eta = f(\p )$, for some $\p \in {}^{<\o }2$. Since $\la \p , \p^\smallfrown k \ra \in \r$, for $k\in\{0,1\}$,
we have $\la f(\p ), f( \p^\smallfrown k )\ra \in \r$ and, hence,
$f( \p^\smallfrown k )=f( \p)^\smallfrown j_k =f(\emptyset )^\smallfrown \eta^\smallfrown j_k$, where $j_k\in \{ 0,1 \}$.
Since $f$ is an injection we have $j_0\neq j_1$ and, hence, $f(\emptyset )^\smallfrown \eta^\smallfrown 0 $ and $f(\emptyset )^\smallfrown \eta^\smallfrown 1$
are elements of $A$. So $A=A_{f(\emptyset )}$ and the proof of (\ref{EQ4071}) is finished.

Using  (\ref{EQ4071}) it is easy to see that $\la {}^{<\o }2 , \supset \ra \cong _F\la \P (G_{{}^{<\o }2} ),\subset \ra $,
where $F(\f )=A_\f$.
\end{ex}
\section{Indivisible structures. Forcing with quotients}
A relational structure $\X =\la X, \{ \r _i :i\in I \} \ra$ is called {\it indivisible} iff for each partition
$X=A \cup B$ we have $\X \hookrightarrow A$ or $\X\hookrightarrow B$.
The aim of this section is to locate indivisible structures in our diagram.
\begin{te}\rm \label{T4072}
A relational structure $\X $ is indivisible iff $\I _{\X }$ is an ideal in $P(X)$.
\end{te}
\dok
Let $\X$ be a indivisible structure. Clearly, $\emptyset \in \I _\X \not\ni X$ and
$I' \subset I \in \I _\X$ implies $I' \in \I _\X$. Suppose that $I\cup J \not\in \I _\X $, for some
$I,J\in \I _\X$. Then $C\subset I\cup J$, for some $C\in \P (\X)$ and $C=(C\cap I)\cup (C \cap (J\setminus I))$.
Since $C\cong \X$, $C$ is indivisible and, hence, there is $A\in \P (C)\subset \P (\X )$ such that
$A\subset C\cap I$ or $A\subset C\cap (J\setminus I)$, which is impossible because $I,J\in \I _\X$.
Thus $\I _{\X }$ is an ideal.

Let $\X$ be a divisible and let $X=A\cup B$ be a partition such that $\X \not\hookrightarrow A$ and $\X\not\hookrightarrow B$.
Then $A,B\in \I _\X$ and, clearly, $A\cup B\not \in \I _\X$. Thus $\I _{\X }$ is not an ideal.
\hfill $\Box$
\begin{te}\rm \label{T4073}
If $\X =\la X , \{ \r_i :i\in I \}\ra$ is an indivisible relational structure, then

(a) $\sm \la \P (\X ), \subset \ra = \la \P (\X ), \subset _{\I _\X} \ra $, where
$A \subset _{\I _\X} B \Leftrightarrow A\setminus B \in \I _\X $;

(b) $\sq \la \P (\X ), \subset \ra $ is isomorphic to a dense subset of  $ \la (P (X )/\!\! =_{\I _\X})^+, \leq _{\I _\X} \ra $. Hence
the poset $\la \P (\X ) , \subset \ra$ is forcing equivalent to $(P(X)/\I _{\X })^+$.

\end{te}
\dok
(a)
Let $A\setminus B\in \I_\X$. If $C \in \P (\X )$ and $C\subset A$, then $C\setminus B\in \I _\X$ and, since $\I _\X $ is an ideal and $C\not\in \I _\X$,
we have $C\cap B\not\in \I _\X$ and, hence, $D\subset C\cap B$, for some $D\in \P (\X)$. By  (\ref{EQ4002}) we have $ A\leq ^* B$.

If $A\setminus B\not\in \I_\X$, then $C\subset A\setminus B$, for some $C\in \P (\X )$ and $C\cap B=\emptyset$ so,
by (\ref{EQ4002}), we have $\neg A\leq ^* B$.

(b) By (a) and the definition of the separative quotient, we have  $\sq \la \P (\X ), \subset \ra = \la \P (\X) /{=^*}, \trianglelefteq \ra$, where for $A,B\in \P (\X )$,
\begin{equation}\label{EQ4027}
A=^* B \Leftrightarrow A\bigtriangleup B \in \I_\X
\;\;\mbox{ and }  \;\;
[A]_{=^*}\trianglelefteq [B]_{=^*}\Leftrightarrow A\setminus B \in \I_\X .
\end{equation}
We show that $\la \P (\X) /{=^*}, \trianglelefteq \ra \hookrightarrow _f  \la (P(X)/\I _\X )^+ ,\leq _{\I _\X }\ra$, where $f([A]_{=^*})=[A]_{=_{\I _\X}}$.
By  (\ref{EQ4027}) and (a),
$[A]_{=^*}= [B]_{=^*}$ iff
$A=^* B$ iff $A \bigtriangleup B \in \I _\X$ iff $A=_{\I _\X} B$ iff $[A]_{=_{\I _\X}}=[B]_{=_{\I _\X}}$
iff $f([A]_{=^*})=f([B]_{=^*})$ and $f$ is a well defined injection.

$f$ is a strong homomorphism since
$[A]_{=^*}\trianglelefteq [B]_{=^*}$ iff
$A\setminus B \in \I _\X$ iff $[A]_{=_{\I _\X}} \leq _{\I _\X } [B]_{=_{\I _\X}}$
iff $f([A]_{=^*})\leq _{\I _\X }  f([B]_{=^*})$.

We prove that $f[\P (\X) /{=^*}]$ is a dense subset of $(P (X )/\!\! =_{\I _\X})^+$.
If $[S]_{=_{\I _\X}}\in (P (X )/\!\! =_{\I _\X})^+$, then $S\not\in \I _\X$ and there is $A\in \P (\X )$ such that
$A\subset S$. Hence $A \subset _{\I _\X} S$ and $f([A]_{=^*})=[A]_{=_{\I _\X}} \leq _{\I _\X } [S]_{=_{\I _\X}}$.

By Fact \ref{T4055}(a) these three posets are forcing equivalent.
\kdok

\noindent
Confirming a conjecture of Fra\"{\i}ss\'{e} Pouzet proved that each countable indivisible structure
contains two disjoint copies of itself \cite{Pouz}. This is, essentially, the statement (a) of the following theorem but, for completeness, we include a proof.

\begin{te}\rm \label{T4074}
If $\X =\la \o , \{ \r _i :i\in I \} \ra$ is a countable indivisible structure, then

(a) $\la \P (\X ) , \subset \ra$ is an atomless partial order (Pouzet);

(b) $|\P (\X )| ={\mathfrak c}$;

(c) $|\sq \la \P (\X ), \subset \ra |>\o$.
\end{te}
\dok
(a) Suppose that $\la \P (\X ) , \subset \ra$ is not atomless. Then, by Theorem \ref{T4013}(d),
$\U =\P (\X )\up $ is a uniform filter on
$\o$. Since $\X$ is indivisible, for each $A\subset \o$ there is $C\in \P (\X )$ such that $C\subset A$ and, hence,
$A\in \U$, or $C\subset \o\setminus A$, and, hence, $\o \setminus A \in \U$. Thus $\P (\X )\up$ is a uniform ultrafilter on $\o$ and,
by a well known theorem of Sierpi\'nski, does not have the Baire property (see e.g. \cite{Kech}, p.\ 56). A contradiction to Theorem \ref{T4052}.

(b) Suppose that $|\P (\X )|<{\mathfrak c}$. Then, by (a) and Theorem \ref{T4052}, we have  $|\P (\X )|=\o$ and, hence,
$\P (\X ) = \{ C_n :n\in \o \}\subset [\o ]^\o$. Since each countable subfamily of $[\o ]^\o$ can be reaped,
there is $A\in [\o ]^\o$ such that $|C_n \cap A|= |C_n \setminus A|=\o$, for each $n\in \o$, and, hence, neither $A$ nor $\o \setminus A$
contain an element of $\P (\X)$, which contradicts the assumption that $\X$ is indivisible.

(c) This is Theorem 3.12 of \cite{Kursq}.
\hfill $\Box$
\section{Embedding-maximal structures}
A relational structure $\X$ will be called {\it embedding-maximal} iff $\P (\X)=[X]^{|X|}$. In this section we characterize
countable embedding-maximal structures and obtain more information on the structures which do not have this property.
If $\P=\la P ,\leq \ra$ is a partial order, a set $S\subset P$ is {\it somewhere dense} in $\P$ iff
there is $p\in P $ such that for each $q\leq p$ there is $s\in S$ satisfying $s\leq q$.
Otherwise, $S$ is {\it nowhere dense}.
\begin{te}\rm \label{T4075}
For a countable binary relational structure $\X=\la \o , \r \ra$ the following conditions are equivalent:

\vspace{-3mm}

\begin{itemize}\itemsep=-2mm
\item[(a)] $\P (\X )=[\o ]^\o$;
\item[(b)] $\P (\X )$ is a dense set in $\la [\o ]^\o , \subset \ra $;
\item[(c)] $\X=\la \o , \r \ra$ is isomorphic to one of the following relational structures:

\vspace{-3mm}

\begin{itemize}\itemsep=-1mm
\item[1] The empty relation, $\la \o , \emptyset \ra$,
\item[2] The complete graph, $\la \o , \o ^2 \setminus \Delta _\o \ra$,
\item[3] The natural strict linear order on $\o$, $\la \o , < \ra$,
\item[4] The inverse of the natural strict linear order on $\o$, $\la \o , < ^{-1}\ra$,
\item[5] The diagonal relation, $\la \o , \Delta _\o \ra$,
\item[6] The full relation, $\la \o , \o ^2 \ra$,
\item[7] The natural linear order on $\o$, $\la \o , \leq \ra$,
\item[8] The inverse of the natural linear order on $\o$, $\la \o , \leq ^{-1}\ra$;
\end{itemize}

\vspace{-2mm}

\item[(d)] $\P (\X )$ is a somewhere dense set in $\la [\o ]^\o , \subset \ra $;
\item[(e)] $\I _\X =\Fin$.
\end{itemize}

\vspace{-3mm}

\noindent
Then the poset $\sq \la \P (\X ), \subset \ra = (P(\o)/\Fin )^+$ is atomless and $\s$-closed.
\end{te}
\dok
The implication (a) $\Rightarrow$ (b) is trivial and it is easy to check (c) $\Rightarrow$ (a).

(b) $\Rightarrow$ (c). Let $\P (\X )$ be a dense set in $\la [\o ]^\o , \subset \ra $.

\vspace{2mm}
\noindent
{\it Claim 1.}
The relation $\r$ is reflexive or irreflexive.

\vspace{2mm}
\noindent
{\it Proof of Claim 1.}
If $R=\{ x\in \o : x\r x\}\in [\o ]^\o$, then there is $C\subset R$ such that $\la \o ,\r  \ra \cong \la C,\r _C\ra$
and, since $\r _C$ is reflexive, $\r$ is reflexive as well. Otherwise we have $I=\{ x\in \o : \neg x\r x\}\in [\o ]^\o$ and, similarly,
$\r$ must be irreflexive.

\vspace{2mm}
\noindent
{\it Claim 2.}
If the relation $\r$ is irreflexive, then the structure $\la \o , \r \ra$ is isomorphic to one of the structures 1 - 4 from (c).

\vspace{2mm}
\noindent
{\it Proof of Claim 2.}
Clearly, $[\o ]^2=K_0 \cup K_1 \cup K_2 \cup K_3$, where the sets

$K_0=\{ \{ x,y \} \in [\o ]^2 : \neg x\r y \land \neg y\r x\}$,

$K_1=\{ \{ x,y \} \in [\o ]^2 :  x\r y \land  y\r x\}$,

$K_2=\{ \{ x,y \} \in [\o ]^2 :  x\r y \land \neg  y\r x \land x<y\}$,

$K_3=\{ \{ x,y \} \in [\o ]^2 :  x\r y \land \neg  y\r x \land x>y\}$,

\noindent
are disjoint. By Ramsey's theorem there are $H\in [\o ]^\o$ and $i\in \{ 0,1,2,3 \}$
such that $[H]^2 \subset K_i$. Since $\P (\X )$ is a dense set in $\la [\o ]^\o , \subset \ra $, there is $C\subset H$
such that
\begin{equation}\label{EQ4003}
\la \o ,\r  \ra \cong \la C,\r _C\ra.
\end{equation}

If $[H]^2 \subset K_0$, then for different $x,y\in C$ we have $\neg x\r y$ and, since $\r$ is irreflexive,
$\r _C=\emptyset$. By (\ref{EQ4003}) we have $\r =\emptyset$.

If $[H]^2 \subset K_1$, then for different $x,y\in C$ we have $x\r y $ and $y \r x$. So, since $\r$ is irreflexive,
$\r _C=C^2 \setminus \Delta _C$, that is  the structure $\la C,\r _C\ra$ is a countable complete graph. By (\ref{EQ4003}) we have $\r =\o ^2 \setminus \Delta _\o$.

If $[H]^2 \subset K_2$, then for different $x,y\in C$ we have
\begin{equation}\label{EQ4004}
(x\r y \land \neg y\r x \land x< y) \lor (y\r x \land \neg x\r y \land y< x).
\end{equation}
Let us prove that for each $x,y\in C$
\begin{equation}\label{EQ4005}
x\r y \Leftrightarrow x<y.
\end{equation}
If $x=y$, then, since $\r$ is irreflexive, we have $\neg x\r y$ and, since $\neg x<y$, (\ref{EQ4005}) is true.

If $x<y$, by (\ref{EQ4004}) we have $x\r y$ and (\ref{EQ4005}) is true.

If $x>y$, by (\ref{EQ4004}) we have $\neg x\r y$ and, since $\neg x<y$, (\ref{EQ4005}) is true again.

\noindent
Since (\ref{EQ4005}) holds for each $x,y\in C$ we have $\r _C = <_C$. Clearly $\la C, <_C\ra \cong \la \o , < \ra$, which,
together with (\ref{EQ4003}), implies $\la \o ,\r  \ra \cong\la \o , < \ra$.

If $[H]^2 \subset K_3$, then as in the previous case we show that $\la \o ,\r  \ra \cong\la \o , <^{-1} \ra$.

\vspace{2mm}
\noindent
{\it Claim 3.}
If the relation $\r$ is reflexive and $\Y =\la \o , \r \setminus \Delta _\o \ra$, then

(i) $\P (\Y )$ is a dense set in $\la [\o ]^\o , \subset \ra $;

(ii) The structure $\la \o , \r \ra$ is isomorphic to one of the structures 5 - 8 from (c).

\vspace{2mm}
\noindent
{\it Proof of Claim 3.}
(i) Let $A\in [\o ]^\o$, $C\subset A$ and $\la \o ,\r \ra \cong _f \la C, \r _c \ra$. Then, since $f$ is an isomorphism, we have
$\la x_1 , x_2 \ra \in \r \setminus \Delta _\o$ iff $\la x_1 , x_2 \ra \in \r \land x_1 \neq x_2$ iff
$\la f(x_1) , f(x_2 )\ra \in \r _C \land f( x_1 ) \neq f(x_2 )$ iff $\la f(x_1) , f(x_2 )\ra \in \r _C \setminus \Delta _\o =(\r \setminus \Delta _\o)_C$.
Thus $\la \o ,\r \setminus \Delta _\o \ra \cong _f \la C, (\r\setminus \Delta _\o)_C \ra$, which implies $C\in \P (\Y )$.

(ii) Since $\r \setminus \Delta _\o$ is an irreflexive relation, by (i) and Claim 2 the structure $\la \o ,\r \setminus \Delta _\o \ra$ is isomorphic
to one of the structures 1 - 4. Hence the structure $\la \o ,\r  \ra$ is isomorphic
to one of the structures 5 - 8.

(b) $\Leftrightarrow$ (e). Since $\I _\X = P(\o )\setminus (\P (\X )\up)$ we have:
$\P (\X )$ is a dense set in $\la [\o ]^\o , \subset \ra $ iff
$\P (\X )\up = [\o ]^\o$ iff $\I _\X =\Fin$.

(b) $\Rightarrow$ (d) is trivial.

(d) $\Rightarrow$ (b) Let $\P (\X )$ be dense below $A\in [\o ]^\o $. Then there are $C\subset A$ and $f$
such that $\X \cong _f \la C, \r _C \ra$ and, by the assumption,
\begin{equation}\label{EQ4006}
\forall B \in [C]^\o \;\; \exists D\in \P (\X )\;\; D\subset B .
\end{equation}
For $S\in [\o ]^\o$ we have $f[S]\in [C]^\o$ and, by (\ref{EQ4006}), there is $D\subset f[S]$ such that $\X \cong \la D, \r _D\ra$.
Since $f$ is an injection we have $f^{-1}[D]\subset S$;  $D\subset f[S]$ implies $f[f^{-1}[D]]=D$ and, since $f$ is an isomorphism, $\la f^{-1}[D], \r _{f^{-1}[D]}\ra \cong _{f | f^{-1}[D]} \la D, \r _D\ra$ and, hence, $f^{-1}[D]\in \P (\X )$. Thus
$\P (\X )$ is a dense set in $\la [\o ]^\o , \subset \ra $.
\hfill $\Box$

\begin{cor}\rm \label{T4058}
If $\X=\la \o , \r \ra$ is a countable binary relational structure, then

(a) $\P (\X )=[\o ]^\o$ or $\P (\X )$ is a nowhere dense set in $\la [\o ]^\o , \subset \ra $;

(b) If $\X$ is indivisible, then $\I _\X=\Fin$ or $\I _\X$ is a tall ideal (that is, for each $S\in [\o ]^\o$ there is $I\in \I _\X\cap [S]^\o$).
\end{cor}
\dok
(b) If  $\I _\X\neq \Fin$, then, by Theorem \ref{T4075}, $\P (\X )$ is a nowhere dense subset of $[\o ]^\o$, so for $S\in [\o ]^\o$
there is $I\in [S]^\o$ such that $A\subset I$, for no $A\in \P (\X )$, which means that $I\in \I _\X$.
\hfill $\Box$
\section{Embeddings of disconnected structures}
If $\X _i=\la X_i, \r _i \ra$, $i\in I$, are binary relational structures  and $X_i \cap X_j =\emptyset$, for
different $i,j\in I$, then the structure $\bigcup _{i\in I} \X _i =\la \bigcup _{i\in I} X_i , \bigcup _{i\in I} \r _i\ra$ will be called
the {\it disjoint union} of the structures $\X _i$, $i\in I$.

If $\la X,\r \ra$ is a binary structure, then
the transitive closure $\r _{rst}$ of the relation $\r _{rs} =\Delta _X \cup \r \cup \r ^{-1}$ (given by $x \;\r_{rst} \;y$ iff there are $n\in \N$ and
$z_0 =x , z_1, \dots ,z_n =y$ such that $z_i \;\r _{rs} \;z_{i+1}$, for each $i<n$)
is the minimal equivalence relation on $X$ containing $\r$.
In the sequel the relation $\r _{rst}$ will be denoted by $\sim _\r$ or $\sim$. Then for $x\in X$
the corresponding element of the quotient $X/\!\!\sim$ will be denoted by $[x]_{\sim _\r}$ or $[x]_\sim$ or only by $[x]$, if the context admits,  and called the
{\it component} of $\la X,\r \ra$ containing $x$. The structure $\la X,\r \ra$ will be called {\it connected}
iff $|X/\!\!\sim |=1$. The main result of this section is Theorem \ref{T4015} describing embeddings of disconnected structures and providing several constructions
in the sequel.
\begin{lem}\rm\label{T4001}
Let $\la X,\r \ra =\la \bigcup _{i\in I} X_i , \bigcup _{i\in I} \r _i\ra$ be a disjoint union of binary structures. Then for each $i\in I$ and each $x\in X_i$ we have

(a) $[x]\subset X_i$;

(b) $[x]=X_i$, if $\la X_i ,\r _i \ra$ is a connected structure.
\end{lem}
\dok
(a) Let   $y\in [x]$ and
$z_0 =x , z_1, \dots ,z_n =y\in X$, where $z_k \;\r _{rs} \;z_{k+1}$, for each $k<n$.
Using induction we show that $z_k \in X_i$, for each $k\leq n$. Suppose that $z_k \in X_i$. Then $z_k \;\r _{rs} \;z_{k+1}$ and,
if $z_k = z_{k+1}$, we are done.
If $\la z_k, z_{k+1}\ra \in \r $, there is $j\in I$ such that
$\la z_k, z_{k+1}\ra \in \r _j \subset X_j\times \X_j$ and, since $z_k \in X_i$, we have $j=i$ and, hence,
$z_{k+1} \in X_i$.
If $\la z_k, z_{k+1}\ra \in \r ^{-1} $, then $\la z_{k+1}, z_k\ra \in \r $ and, similarly,  $z_{k+1} \in X_i$ again.

(b) Let $\la X_i ,\r _i \ra$ be a connected structure and $y\in X_i$.
Then $x \sim _{\r _i }y$ and, hence, there are $z_0 =x , z_1, \dots ,z_n =y\in X_i$, where for each $k<n$ we have $z_k \;(\r _i) _{rs} \;z_{k+1}$, that is
$z_k = z_{k+1} \lor z_k \;\r _i\; z_{k+1}\lor z_k \;(\r _i)^{-1}\; z_{k+1}$, which implies $z_k \;\r  _{rs} \;z_{k+1}$.
Thus $y\sim _\r x$ and, hence, $y\in [x]$.
\hfill $\Box$
\begin{prop}\rm\label{T4002}
If $\la X,\r \ra$ is a binary structure, then
$ \la \bigcup _{x\in X}[x],\bigcup _{x\in X}\r_{[x]} \ra$
is the unique representation of $\la X,\r \ra$ as a disjoint union of connected relations.
\end{prop}
\dok
Clearly $X=\bigcup _{x\in X}[x]$ is a partition of $X$ and $\bigcup _{x\in X}\r_{[x]} \subset \r$.
If $\la x,y\ra \in \r$, then $x\sim y$, which implies $x,y\in [x]$. Hence $\la x,y\ra \in \r \cap ([x] \times [x])=\r_{[x]}$ and we have
$\r = \bigcup _{x\in X}\r_{[x]} $.

We show that the structures $\la [x],\r_{[x]} \ra$, $x\in X$, are connected. Let $y\in [x]$
and $z_0 =x , z_1, \dots ,z_n =y\in X$, where $z_k \;\r _{rs} \;z_{k+1}$, for each $k<n$.
Using induction we show that
\begin{equation}\label{EQ4007}
\forall k \leq n \;\; z_k \in [x].
\end{equation}
Suppose that $z_k \in [x]$. Then $z_k \;\r _{rs} \;z_{k+1}$ and,
if $z_k = z_{k+1}$, we are done.
If $\la z_k, z_{k+1}\ra \in \r $, there is $u\in X$ such that
$\la z_k, z_{k+1}\ra \in \r _{[u]} \subset [u]\times [u]$ and, since $z_k \in [x]$, we have $[u]=[x]$ and, hence,
$z_{k+1} \in [x]$.
If $\la z_k, z_{k+1}\ra \in \r ^{-1} $, then $\la z_{k+1}, z_k\ra \in \r $ and, similarly,  $z_{k+1} \in [x]$ again.

For each $k<n$ we have $\la z_k, z_{k+1}\ra \in \Delta _X \cup \r \cup \r ^{-1}$ so, by (\ref{EQ4007}),
$\la z_k, z_{k+1}\ra \in \Delta _{[x]} \cup \r _{[x]} \cup \r _{[x]} ^{-1} = (\r _{[x]})_{rs}$.
Thus $x\sim _{\r _{[x]}}y$ and, since the relation $\sim _{\r _{[x]}}$ is symmetric, $y\sim _{\r _{[x]}}x$, for each $y\in [x]$.
Since the relation $\sim _{\r _{[x]}}$ is transitive, for each $y,z\in [x]$ we have $y\sim _{\r _{[x]}}z$ and, hence, $\la [x],\r_{[x]} \ra$ is a connected
structure.

For a proof of the uniqueness of the representation, suppose that
$\la X,\r \ra =\la \bigcup _{i\in I} X_i , \bigcup _{i\in I} \r _i\ra$ is a disjoint union, where the structures $\la X_i ,\r _i \ra$, $i\in I$, are
connected. By Lemma \ref{T4001}(b), for $i\in I$ and  $x\in X_i$ we have $X_i=[x]$ and, hence, $\r _i = \r \cap (X_i \times X_i) =\r \cap ([x]\times [x])=\r _{[x]}$.
Thus $\la X_i , \r _i \ra = \la  [x], \r _{[x]}\ra$.
On the other hand, if $x\in X$, then $x\in X_i$, for some $i\in I$, and, similarly,  $\la  [x], \r _{[x]}\ra =\la X_i , \r _i \ra$.
Consequently we have
$\{ \la X_i , \r _i \ra :i\in I \} =\{ \la  [x], \r _{[x]}\ra :x\in X \}$.
\hfill $\Box$
\begin{prop}\rm\label{T4010}
Let $\la X,\r \ra$ be a binary relational structure and $\r ^c = (X\times X)\setminus \r$ the complement of $\r$. Then

(a) At least one of the structures $\la X,\r \ra$ and $\la X,\r ^c \ra$ is connected;

(b) $\Emb \la X,\r \ra = \Emb \la X,\r^c \ra$ and $\P \la X,\r \ra = \P \la X,\r^c \ra$.
\end{prop}
\dok
(a)
Suppose that the structure $\X=\la X,\r \ra$ is disconnected. Then, by Proposition \ref{T4002},
$\X$ is the disjoint union of connected structures $\X _i=\la X_i ,\r_i \ra$, $i\in I$, and we show that  $\la X,\r ^c \ra$ is connected.
Let $x,y\in X$.
If $x\in X_i$ and $y\in X_j$, where $i\neq j$, then $x \not\sim _{\r }y$, which implies $\la x,y \ra \not\in \r$, thus
$\la x,y \ra \in \r ^c$ and, hence, $x \sim _{\r ^c} y$.
Otherwise, if $x,y\in X_i$, for some $i\in I$, then we pick $j\in I\setminus \{ i\}$
and $z\in X_j$ and, as in the previous case,  $x \sim _{\r ^c} z$ and $y \sim _{\r ^c} z$ and, since $\sim _{\r ^c} $ is an equivalence relation,
$x \sim _{\r ^c} y$ again.

(b) If $f\in \Emb \la X,\r \ra$, then  $f$ is an injection and for each $x,y\in X$ we have $\la x,y \ra \in \r \Leftrightarrow \la f(x),f(y )\ra \in \r$,
that is $\la x,y \ra \in \r ^c\Leftrightarrow \la f(x),f(y )\ra \in \r ^c$ and, hence, $f\in \Emb \la X,\r ^c\ra$.
The another implication has a similar proof. Now
$\P \la X,\r \ra
= \{ f[X] :f\in \Emb \la X,\r \ra \}
= \{ f[X] :f\in \Emb \la X,\r ^c\ra  \}
=\P \la X,\r^c \ra$.
\hfill $\Box$
\begin{lem}\rm\label{T4005}
Let $\la X, \r \ra$ and $\la Y, \t \ra$ be binary structures and $f:X \rightarrow Y$ an embedding. Then
for each $x_1 , x_2 , x \in X$

(a) $x_1 \r _{rs} x_2 \Leftrightarrow f(x_1 ) \t _{rs} f( x_2 )$;

(b) $x_1 \sim _\r  x_2 \Rightarrow f(x_1 ) \sim _\t f( x_2 )$;

(c) $f[[x]] \subset [f(x)]$;

(d) $f\,|\;[x]:[x]\rightarrow f[[x]]$ is an isomorphism.

\noindent
If, in addition, $f$ is an isomorphism, then

(e) $x_1 \sim _\r  x_2 \Leftrightarrow f(x_1 ) \sim _\t f( x_2 )$;

(f) $f[[x]]= [f(x)]$;

(g) $\la X, \r \ra$ is connected iff  $\la Y, \t \ra$ is connected.
\end{lem}
\dok
(a) Since $f$ is an injection and a strong homomorphism we have $x_1 \; \r _{rs} x_2$ iff
$x_1=x_2 \lor x_1 \; \r \; x_2 \lor x_2 \; \r \; x_1$ iff $f(x_1)=f(x_2) \lor f(x_1) \; \r \; f(x_2) \lor f(x_2) \; \r \; f(x_1)$
iff $f(x_1) \; \t _{rs} f(x_2)$.

(b) If $x_1 \sim _\r  x_2 $, then there are $z_0, z_1,\dots , z_n \in X$
such that
$x_1=z_0 \;\r _{rs }\; z_1 \;\r _{rs }$ $ \dots  \r _{rs }\; z_n =x_2$
and, by (a),
$f(x_1)=f(z_0) \;\t _{rs }\; f(z_1) \;\t _{rs } \dots  \t _{rs }\; f(z_n) =f(x_2)$ and,
hence, $f(x_1 ) \sim _\t f( x_2 )$.

(c) If $x' \in [x]$, then $x' \sim _{\r } x$ and, by (b), $f(x') \sim _{\t } f(x)$ so
$f(x')\in [f(x)]$.

(d) Clearly, $f|[x]$ is a bijection. Since $f$ is a strong homomorphism, for $x_1,x_2 \in [x]$ we have $x_1 \; \r  x_2$ iff
$f(x_1) \; \t  f(x_2)$ iff $(f|[x])(x_1) \; \t  (f|[x])(x_2)$.

(e) The implication ``$\Rightarrow $" is proved in (b). If $f(x_1 ) \sim _\t f( x_2 )$, then, applying (b) to $f^{-1}$ we obtain  $x_1 \sim _\r  x_2$.

(f) The inclusion ``$\subset $" is proved in (b). Let $y\in [f(x)]$, that is $y\sim _\t f(x)$. Since $f$ is a bijection there is $x'\in X$
such that $y=f(x')$ and, by (e), $x' \sim _\r x$, that is $x' \in [x]$. Hence $y\in f[[x]]$.

(g) follows from (e).
\hfill $\Box$
\begin{te}\rm\label{T4015}
Let $\X _i = \la X_i , \r _i \ra , i\in I$, and $\Y _j = \la Y_j , \s _j \ra , j\in J$, be two families of disjoint connected
binary structures and $\X$ and $\Y$ their unions.
Then

(a) $F:\X  \hookrightarrow \Y $ iff there are $f:I\rightarrow J$ and $g_i : \X _i \hookrightarrow \Y _{f(i)}$, $i\in I$,
such that $F=\bigcup _{i\in I} g_i$ and
\begin{equation}\label{EQ4080}
\forall \{ i_1 ,i_2 \}\in [I]^2 \; \forall x _{i_1} \in X_{i_1} \; \forall  x _{i_2} \in X_{i_2}
       \; \neg \; g_{i_1}(x_{i_1})\; \s _{rs} \;g_{i_2}(x_{i_2}).
\end{equation}

(b) $C\in \P (\X )$ iff there are $f:I\rightarrow I$ and $g_i : \X _i \hookrightarrow \X _{f(i)}$, $i\in I$,
such that $C= \bigcup _{i\in I} g_i [X_i]$ and
\begin{equation}\label{EQ4009}
\forall \{i,j\} \in [I]^2 \;\; \forall x \in X_i \;\; \forall y \in X_j \;\;\neg \; g_i (x)\; \r _{rs} \;g_j(y).
\end{equation}
\end{te}
\dok
(a)
($\Rightarrow$) Let $F:\X  \hookrightarrow \Y $. By Proposition \ref{T4002}, the sets $X_i$, $i\in I$, are components of $\X$ and
$Y_j$, $j\in I$, are components of $\Y$.
By Lemma \ref{T4005}(c), for $i\in I$ and $x\in X_i$ we have $F[[x]]\subset [F(x)]$ so there is (unique) $f(i)\in J$, such that $F[X_i]\subset Y_{f(i)}$.
By Lemma \ref{T4005}(d), $F|X_i : X_i \rightarrow F[X_i] \subset Y_{f(i)}$ is an isomorphism and, hence, $g_i : \X _i \hookrightarrow \Y _{f(i)}$, where the mapping
$g_i : X_i \rightarrow Y_{f(i)}$ is given by $g_i(x)=F(x)$. Clearly $f:I\rightarrow J$ and  $F=\bigcup _{i\in I} g_i$.
Suppose that $g_{i_1}(x_{i_1})\; \s _{rs} \;g_{i_2}(x_{i_2})$, that is $F(x_{i_1})\; \s _{rs} \;F(x_{i_2})$,
for some different $i_1 , i_2 \in I$ and some $x _{i_1} \in X_{i_1}$ and $ x _{i_2} \in X_{i_2}$. Then, by Lemma \ref{T4005}(a),
$x_{i_1}\; \r _{rs} \;x_{i_2}$ and, hence, $x_{i_1} \sim _\r x_{i_2}$, which is not true, because $x_{i_1}$ and $x_{i_2}$ are elements of different
components of $\X$.

($\Leftarrow$) Let $F=\bigcup _{i\in I} g_i$, where the functions $f:I\rightarrow J$ and $g_i : \X _i \hookrightarrow \Y _{f(i)}$, $i\in I$, satisfy
the given conditions.

Let $u,v\in X$, where $u\neq v$.
If $u,v\in X_i$ for some $i\in I$ then, since $g_i$ is an injection, we have $F(u)=g_i(u)\neq g_i (v)=F(v)$. Otherwise $u\in X_{i_1}$ and $v\in X_{i_2}$, where
$i_1 \neq i_2$ and, by the assumption, $\neg \; g_{i_1}(u)\; \s _{rs} \;g_{i_2}(v)$, which implies $ g_{i_1}(u)\neq g_{i_2}(v)$ that is $F(u)\neq F(v)$. Thus
$F$ is an injection.

In order to prove that $F$ is a strong homomorphism we take $u,v\in X$ and prove
\begin{equation}\label{EQ4008}
u\; \r \; v \Leftrightarrow F(u)\; \s F(v) .
\end{equation}
If $u,v \in X_i$, for some $i\in I$, then we have: $u\; \r \; v$ iff $u\; \r _i \; v$ (since $\r _{X_i}=\r _i$) iff
$g_i (u) \; \s _{f(i)} \; g_i (v)$ (because $g_i : \X _i \hookrightarrow \Y _{f(i)}$) iff $g_i (u) \; \s  \; g_i (v)$
(since $\s _{Y_{f(i)}}=\s _{f(i)}$) iff $F(u)\;\s \; F(v)$ (because $F\upharpoonright X_i =g_i$). So (\ref{EQ4008}) is true.

If $u\in X_{i_1}$ and $v\in X_{i_2}$, where $i_1 \neq i_2$, then $\neg u\; \r \;v$, because $u$ and $v$ are in different components of $X$.
By the assumption we have $\neg \; g_{i_1}(u)\; \s _{rs} \;g_{i_2}(v)$, which implies $\neg \; g_{i_1}(u)\; \s \;g_{i_2}(v)$, that is
$\neg \; F(u)\; \s \;F(v)$. So (\ref{EQ4008}) is true again.

(b) follows from (a) and the fact that $C\in \P (\X )$ iff there is $F:\X \hookrightarrow \X$ such that $C=F[X]$.
\hfill $\Box$
\section{Embedding-incomparable components}
Two structures $\X$ and $\Y$ will be called {\it embedding-incomparable} iff $\X \not\hookrightarrow\Y$ and $\Y \not\hookrightarrow\X$.
We will use the following fact.
\begin{fac}  \rm \label{T4042}
Let $\P , \Q $ and $\P _i$, $i\in I$, be partial orderings. Then

(a) If $\P \cong \Q$, then $\sm \P \cong  \sm \Q $ and $\sq \P \cong  \sq \Q$;

(b) $\sm (\prod _{i\in I}\P _i) = \prod _{i\in I}\sm \P _i$;

(c) $\sq (\prod _{i\in I}\P _i) \cong \prod _{i\in I}\sq \P _i$.
\end{fac}
\begin{te}\rm\label{T4016}
Let $\r $ be a binary relation on a set $X$.
If the components $\X _i =\la X_i , \r _{X_i}\ra $, $i\in I$, of the structure $\X =\la X, \r \ra $ are embedding-in\-com\-pa\-ra\-ble, then

(a) $\la \P (\X ) ,\subset \ra \cong  \prod _{i\in I } \la \P (\X _i), \subset \ra  $;

(b) $\sq \la \P (\X ) , \subset \ra \cong  \prod _{i\in I } \sq \la \P (\X _i) , \subset \ra $.

(c) $\X$ is a divisible structure.
\end{te}
\dok
(a) By Theorem \ref{T4015}(b) and since the structures $\X _i$ are embedding-in\-com\-pa\-ra\-ble,
$C\in \P (\X )$ iff there are embeddings $g_i : \X _i \hookrightarrow \X _i$, $i\in I$,
such that $C= \bigcup _{i\in I} g_i [X_i]$ and $\neg \; g_i (x)\; \r _{rs} \;g_j(y)$,
       for each different $i,j \in I$, each $x \in X_i$ and  each $y \in X_j$.
But, since $i\neq j$,  $x \in X_i$ and $y \in X_j$ implies    $g_i (x)\in X_i$ and  $g_j (y)\in X_j$,
it is impossible that $g_i (x)\; \r _{rs} \;g_j(y)$ and, hence, the last condition is implied by
the condition that  $g_i : \X _i \hookrightarrow \X _i$, for each $i\in I$. Consequently,
$\P (\X )=\{ \bigcup _{i\in I} C_i : \la C_i :i\in I\ra \in \prod _{i\in I } \P (\X _i) \}$ and it is easy to check
that the mapping $f:   \prod _{i\in I } \la \P (\X _i) \subset \ra \rightarrow \la \P (\X ) , \subset \ra$ given by
$f(\la C_i :i\in I\ra )=\bigcup _{i\in I} C_i $ is an isomorphism of posets.

(b) follows from (a) and Fact \ref{T4042}(a) and (c).

(c) The partition $X=X_i \cup (X\setminus X_i)$ witnesses that $\X$ is divisible.
\hfill $\Box$
\section{From $\bf A_1$ to $\bf D_5$}
In this section we show that the diagram on Figure \ref{FIG2} is correct.
The relations between the properties of $X$ and $\P (\X)$ are established in the previous sections.
Since $|\sq \la \P (\X ),\subset \ra|\leq |\P (\X )|$, the classes $B_1$, $C_1$, $D_1$, $C_2$ and $D_2$ are empty
and, since $\sq \la [\o ]^\o , \subset \ra=(P(\o )/\Fin )^+$ is a $\s$-closed atomless poset, the classes
$A_5$, $B_5$ and $C_5$ are empty as well. By Theorem \ref{T4074} we have $A_4=B_4=\emptyset$ and in the sequel we show that
the remaining classes contain some structures. First, the graph $G_{\Z}$ mentioned in the Introduction belongs to $A_1$ and its
restriction to $\N$ to $A_2$. The class $B_2$ contains the digraph constructed in Example \ref{EX4014} and in the following examples we construct some
structures from $A_3$, $B_3$ and $C_3$.
\begin{ex}\rm\label{EX4025}
$\la \P (\X ) ,\subset \ra $ collapses ${\mathfrak c}$ to $\o$ and $\X$ is a divisible structure belonging to $C_3$.
Let $\X =\la X, \r \ra = \la \bigcup _{n\geq 3}G_n',\bigcup _{n\geq 3}\r_n' \ra$, where the sets $G_n'$, $n\geq 3$, are pairwise disjoint and
$\la G_n' , \r _n' \ra\cong \la G_n , \r _n \ra$, where the structure $\la G_n , \r _n \ra$ is the directed graph defined by
$G_n = {}^{<\o }2 \times \{0,1, \dots , n-1\}$ and
\begin{eqnarray*}
\!\!\r_n \!\!& \!\!= \!\!&\!\! \{ \la \la \f , 0\ra , \la \f ^\smallfrown k , 0 \ra\ra : \f \in {}^{<\o }2 \land k\in 2 \}  \cup \\
           &           &\!\! \{ \la \la \f , i\ra , \la \f , j \ra\ra :
           \f \in {}^{<\o }2 \land \la i,j \ra \in \{\la 0,1 \ra , \la 1,2 \ra , \dots , \la n-1 , 0 \ra\}\} .
\end{eqnarray*}
Using the obvious fact that two cycle graphs of different size are embedding incomparable we easily prove that
for different $m,n\geq 3$ the structures  $\la G_m , \r _m \ra$ and $\la G_n , \r _n \ra$ are embedding incomparable as well
so, by (a) of Theorem \ref{T4016},
\begin{equation}\label{EQ4073}\textstyle
\la \P (\X ) ,\subset \ra \cong  \prod _{n\geq 3} \la \P (\la G_n , \r _n \ra ), \subset \ra  .
\end{equation}
Let $n\geq 3$. Like in Example \ref{EX4014} for $\f \in {}^{<\o }2$ let $A_\f =\{\p \in {}^{<\o }2 : \f \subset \p \}$ and
$B_\f =A_\f \times \{ 0,1, \dots , n-1 \}$.
Let us prove that
\begin{equation}\label{EQ4074}
 \P (\la G_n , \r _n \ra  ) =\{ B_\f : \f \in {}^{<\o }2 \} .
\end{equation}
The inclusion ``$\supset$" is evident.
Conversely,
let $B\in \P ( \la G_n , \r _n \ra )$ and $f : \la G_n , \r _n \ra   \hookrightarrow \la G_n , \r _n \ra $, where $B=f[G_n ]$.
Clearly, $\deg (v)\in\{ 4, 5\}$, for each vertex $v\in {}^{<\o }2 \times \{0 \}$, and $\deg (v)=2$, otherwise
Thus, since $f$ preserves degrees of vertices we have $f[{}^{<\o }2 \times \{0 \}]\subset {}^{<\o }2 \times \{0 \}$
and $f\upharpoonright {}^{<\o }2 \times \{0 \}  : {}^{<\o }2 \times \{0 \} \hookrightarrow {}^{<\o }2 \times \{0 \}$.
Since the digraph ${}^{<\o }2 \times \{0 \}$ is isomorphic to the digraph $G_{{}^{<\o }2}$, by
Example \ref{EX4014}, there is $\f \in {}^{<\o }2 $ such that
\begin{equation}\label{EQ4075}
f[{}^{<\o }2 \times \{0 \}]= A_\f \times \{0 \}.
\end{equation}
Now, since each $v\in G_n$ belongs to a unique cycle graph with $n$ vertices and $f$ preserves this property by (\ref{EQ4075}) we have
$B=f[G_n]=B_\f$ and (\ref{EQ4074}) is proved.

By (\ref{EQ4074}), like in  Example \ref{EX4014} we prove that $\la \P (\la G_n , \r _n \ra  ), \subset \ra \cong \la {}^{<\o }2 , \supset \ra$.
Thus, by (\ref{EQ4074}), the poset $\la \P (\X ) ,\subset \ra$  is isomorphic to the direct product $\la {}^{<\o }2 , \supset \ra^\o $
of countably many Cohen posets which collapses ${\mathfrak c}$ to $\o$ (see \cite{Kun}, (E4) on page 294).
The partition $X=G_3 \cup (X\setminus G_3)$ witnesses that $\X$ is a divisible structure.
\end{ex}
\begin{ex}\rm\label{EX4026}
$\la \P (\X ) , \subset \ra$ is an atomic poset of size ${\mathfrak c}$ and $\X \in A_3$.
Let $\X =\la X, \r \ra = \la \bigcup _{n\geq 3}G_n',\bigcup _{n\geq 3}\r_n' \ra$, where the sets $G_n'$, $n\geq 3$, are pairwise disjoint and
$\la G_n' , \r _n' \ra$ is isomorphic to the digraph $\la G_n , \r _n \ra$ given by
$G_n = \o \times \{0,1, \dots , n-1\}$ and
\begin{eqnarray*}
\!\!\r_n \!\!& \!\!= \!\!&\!\! \{ \la \la n , 0\ra , \la n+1 , 0 \ra\ra : n\in \o \}  \cup \\
           &           &\!\! \{ \la \la n , i\ra , \la n , j \ra\ra :
           n\in \o  \land \la i,j \ra \in \{\la 0,1 \ra , \la 1,2 \ra , \dots , \la n-1 , 0 \ra\}\} .
\end{eqnarray*}
As in Example \ref{EX4025} we prove that
for different $m,n\geq 3$ the structures  $\la G_m , \r _m \ra$ and $\la G_n , \r _n \ra$ are embedding incomparable
so, by (a) of Theorem \ref{T4016},
\begin{equation}\label{EQ4076}\textstyle
\la \P (\X ) ,\subset \ra \cong  \prod _{n\geq 3} \la \P (\la G_n , \r _n \ra ), \subset \ra  .
\end{equation}
Let $n\geq 3$. Using the arguments from Example \ref{EX4025} we easily prove that
\begin{equation}\label{EQ4077}
 \P (\la G_n , \r _n \ra  ) =\{ B_k : k \in \o \} ,
\end{equation}
where $B_k =(\o \setminus k ) \times \{ 0,1, \dots , n-1 \}$, for $k\in \o$.

By (\ref{EQ4077}) we have $\la \P (\la G_n , \r _n \ra  ), \subset \ra \cong \la \o , \geq \ra =\o ^*$.
Thus, by (\ref{EQ4076}), the poset $\la \P (\X ) ,\subset \ra$  is isomorphic to the direct product $(\o ^*)^\o $
of countably many copies of $\o ^*$ which is an atomic lattice of size ${\mathfrak c}$.
\end{ex}
\begin{ex}\rm\label{EX4024}
$\sq \la \P (\X ), \subset \ra \cong \la {}^{<\o }2 , \supset \ra$ although $|\P (\X )|={\mathfrak c}$, thus $\X \in B_3$. Let $\Y=\la Y, \r \ra$ be  the digraph considered in
Example \ref{EX4014} and $\Z =\la Z,\s ^c \ra$, where $\la Z, \s \ra$ is isomorphic to the digraph from Example \ref{EX4026} and
$Y\cap Z=\emptyset$. Since $\la Z, \s \ra$ is a disconnected structure, by Proposition \ref{T4010}(a) the structure $\Z$ is connected and, clearly,
$\s ^c = (Z\times Z)\setminus \s$ is a reflexive relation, which implies that the structures $\Y$ and $\Z$ are embedding incomparable.
Thus, by Theorem \ref{T4016}(a), for the structure $\X =\Y \cup \Z$ we have
$\la \P (\X ), \subset \ra \cong  \la \P (\Y ), \subset \ra \times  \la \P (\Z ), \subset \ra$ and, since by Proposition \ref{T4010}(b)
$\P (\Z)=\P (\la Z, \s \ra)$, we have $|\P (\X )|={\mathfrak c}$.

By Theorem \ref{T4016}(b) we have $\sq \la \P (\X ), \subset \ra \cong \sq \la \P (\Y ), \subset \ra \times \sq \la \P (\Z ), \subset \ra$.
Since $\la \P (\Z ) , \subset \ra$ is an atomic poset, by Theorem \ref{T4057}(a) we have  $|\sq \la \P (\Z ) , \subset \ra |=1$ and, hence,
$\sq \la \P (\X ), \subset \ra \cong \la {}^{<\o }2 , \supset \ra \times 1 \cong \la {}^{<\o }2 , \supset \ra$.
\end{ex}
In the sequel we show that the remaining classes are non-empty and give more information about some basic classes of structures.

\vspace{2mm}

{\bf Linear orders.}
A linear order $L$ is {\it scattered} iff it does not contain a dense suborder or, equivalently, a copy of the rationals, $\Q$. Otherwise
$L$ is a {\it non-scattered} linear order. So, if $L$ is a countable linear order, we have the following cases.

{\it Case 1}: $L$ is non-scattered. By \cite{KurTod}, for each non-scattered linear order $L$ the poset $\la \P (L), \subset \ra$ is forcing equivalent to
the two-step iteration $\S \ast \pi$, where $\S$ is the Sacks forcing and $1_\S \Vdash `` \pi $ is a $\sigma$-closed forcing".
If the equality sh$(\S )=\aleph _1$
or PFA holds in the ground model, then the second iterand is forcing equivalent to the poset $(P(\o )/\Fin )^+$ of the
Sacks extension. So, if  $L$ is a countable non-scattered linear order, then forcing by $\la \P (L), \subset \ra$ produces reals.
In addition, $L$ is indivisible. Namely, if $Q$ is a copy of $\Q$ in $L$ and $L=A_0\dot{\cup }A_1$, then,
since $\Q$ is indivisible, there is $k\in \{ 0,1\}$ such that $Q\cap A_k$ contains a copy of $\Q$ and, by the universality of $\Q$, $Q\cap A_k$ contains a copy of $L$ as well. Hence, $L\in C_4$.

{\it Case 2}: $L$ is scattered. By \cite{Kurscatt}
for each countable scattered linear order $L$ the partial ordering
$\sq \la \P (L), \subset \ra$ is atomless and $\sigma$-closed. In particular,
if  $\a $ is a countable ordinal and
$\a =\o ^{\g _n +r_n }s_n + \dots + \o ^{ \g _0 +r_0 }s_0 + k$
its representation in the Cantor normal form, where $k\in \o$, $r_i \in \o$, $s_i \in \N$, $\g _i \in \Lim \cup \{ 1 \}$ and
$\g _n +r_n > \dots > \g _0 +r_0$,
then by \cite{Kurord}
\begin{equation}\label{EQ4078}\textstyle
\sq \la \P (\a ), \subset \ra \cong \prod _{i=0}^n \Big( \Big( \rp ^{r_i}\big( P(\o ^{\g _i} )/ \I _{\o ^{\g _i} }\big)\Big)^+ \Big)^{s_i},
\end{equation}
where, for an ordinal $\b$,
$\I _\b =\{ C\subset \b : \b \not\hookrightarrow C\} $
and, for a poset $\P$, $\rp (\P )$ denotes the reduced
power $\P ^\o / \equiv _{\Fin }$ and $\rp ^{k+1}(\P )= \rp (\rp ^k (\P ))$.
In particular, for $\o \leq \a <\o ^\o$ we have
\begin{equation}\label{EQ4079}\textstyle
\sq \Big( \P \big(\sum _{i=n}^0 \o ^{1+r_i} s_i \big), \subset \Big)
\cong\prod _{i=0}^n \Big(\Big( \rp ^{r_i}\big( P(\o )/ \Fin \big)\Big)^+ \Big)^{s_i}.
\end{equation}
Thus if $L$ is a scattered linear order, then $L\in D_3 \cup D_4 \cup D_5$ and, for example, $\o +\o \in D_3$,
$\o \cdot \o \in D_4 $ and $\o \in D_5$, since an ordinal $\a<\omega _1$ is an indivisible structure iff $\a =\o ^{\b}$, for some ordinal $\b>0$.\\[2mm]
\noindent
So, under the CH, for a countable linear order $L$ the poset $\la \P (L), \subset \ra$ is forcing equivalent to
$\S \ast \pi $, where $1_\S \Vdash `` \pi =(P(\check{\o })/\Fin )^+"$, if $L$ is non-scattered;
and to $(P(\o )/\Fin )^+$, if $L$ is scattered.
But it is consistent that the poset $\la \P (\o + \o ), \subset \ra$ is not forcing equivalent to  $(P(\o )/\Fin )^+$:
by (\ref{EQ4079}) we have $\sq \la \P (\o +\o ), \subset \ra \cong (P(\o )/\Fin )^+ \times (P(\o )/\Fin )^+$ and,
by a result of Shelah and Spinas \cite{SheSpi}, it is consistent that $(P(\o )/\Fin )^+$ and its square are not forcing equivalent.

\vspace{2mm}

{\bf Equivalence relations and similar structures.}
By a more general theorem from \cite{Kurstr} we have: If  $\X _i=\la X_i , \r _{X_i}\ra $, $i\in I$,
are the components of a countable binary structure $\X=\la X, \r\ra$, which is

- either an equivalence relation,

- or a disjoint union of complete graphs,

- or a disjoint union of ordinals $\leq \o$,

\noindent
then
$\sq \la \P (\X ), \subset \ra$ is a $\s$-closed atomless poset.
More precisely, if
$N=\{ |X_i |: i\in I \} $, $N_{\mathrm{fin}}=N\setminus \{ \o \}$, $I_\kappa =\{ i\in I : |X_i |=\kappa \}$, $\kappa \in N$,
and $|I_\o|=\mu$, then the following table describes a forcing equivalent and some cardinal invariants of $\la \P (\X ), \subset\ra$
\begin{center}
{\tiny
\begin{tabular}{c|c|c|c}
$\X$                      & $\sq \la \P (\X ), \subset \ra$ is & $\sq \la \P (\X ), \subset \ra$ is & ZFC $\vdash \sq \la \P (\X ), \subset \ra$   \\
                          &  forcing equivalent to        &                                    &   is ${\mathfrak h}$-distributive            \\[2mm] \hline
                          &                                    &                                    &                                              \\
$N\in [\N ]^{<\o }$ or $|I|=1$              & $(P(\o )/\Fin )^+ $                          & ${\mathfrak t}$-closed              & YES             \\[2mm]
$0<|N_{\mathrm{fin}}|,|I_\o |<\o $            & $((P(\o )/\Fin )^+)^n$                       & ${\mathfrak t}$-closed              & NO              \\[2mm]
$|I_\o |<\o = |N_{\mathrm{fin}}|$&  $(P(\Delta )/\ED )^+  \times  ((P(\o )/\Fin)^+)^\mu$ &  $\sigma$-closed     & NO                     \\[2mm]
$|I_\o|=\o$                                 & $(P(\o \times \o)/(\Fin \times \Fin ))^+$    & $\sigma$-closed, not $\o _2$-closed & NO
\end{tabular}
}
\end{center}

\vspace{2mm}

\noindent
where $\Delta =\{ \la m,n\ra \in \N \times \N : n\leq m\}$ and the ideal $\ED$ in $P(\Delta)$ is defined by
$
\ED = \{ S\subset \Delta : \exists r\in \N \;\; \forall m \in \N \;\; |S\cap (\{ m \} \times \{ 1,2,\dots ,m\})|\leq r  \}.
$

The structure $\X$ is indivisible iff $N\in [\N ]^\o$  or $N=\{ 1 \}$ or $|I|=1$ or $|I_\o |=\o$.

\begin{figure}[h]\label{F4002}
\begin{center}

\unitlength  0.8mm 
\linethickness{0.4pt}
\ifx\plotpoint\undefined\newsavebox{\plotpoint}\fi 


\begin{picture}(85,102)(0,0)


\put(5,5){\line(1,0){45}}
\put(10,15){\line(1,0){70}}
\put(20,20){\line(1,0){5}}
\put(70,20){\line(1,0){5}}
\put(10,40){\line(1,0){70}}
\put(45,60){\line(1,0){30}}
\put(10,65){\line(1,0){70}}
\put(20,85){\line(1,0){25}}
\put(10,90){\line(1,0){70}}
\put(5,100){\line(1,0){3}}
\put(47,100){\line(1,0){3}}
\put(5,5){\line(0,1){95}}
\put(10,15){\line(0,1){75}}
\put(20,20){\line(0,1){65}}
\put(45,60){\line(0,1){25}}
\put(50,25){\line(0,1){75}}
\put(50,5){\line(0,1){12}}
\put(75,20){\line(0,1){40}}
\put(80,15){\line(0,1){75}}



\scriptsize

\put(27,100){\makebox(0,0)[cc]{$\X$ ultrahomogeneous}}
\put(48,20){\makebox(0,0)[cc]{$\X $ equivalence relation}}
\put(32,80){\makebox(0,0)[cc]{$\bigcup _ 1F_\o$}}
\put(32,70){\makebox(0,0)[cc]{$\bigcup _\o F_1$}}
\put(32,53){\makebox(0,0)[cc]{$\bigcup _\o F_\o$}}
\put(62,53){\makebox(0,0)[cc]{$\bigcup _{n\in \o } F_n$}}
\put(32,35){\makebox(0,0)[cc]{$\bigcup _\o  F_2$}}
\put(32,27){\makebox(0,0)[cc]{$\bigcup _2  F_\o$}}
\put(62,30){\makebox(0,0)[cc]{$F_3 \cup \bigcup _\o  F_2$}}
\put(15,30){\makebox(0,0)[cc]{$D_3$}}
\put(15,53){\makebox(0,0)[cc]{$D_4$}}
\put(15,75){\makebox(0,0)[cc]{$D_5$}}
\end{picture}
\end{center}

\vspace{-5mm}

\caption{Equivalence relations on countable sets}
\end{figure}
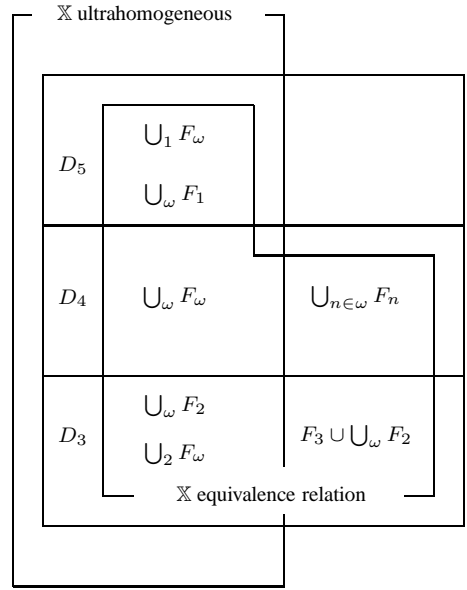

Thus if $\X$ is a countable equivalence relation, then $\X\in D_3 \cup D_4 \cup D_5$ and some examples of such structures are given in the
diagram in Figure \ref{F4002}. We remark that, if $F_\k$ denotes the full relation on a set of size $\k$, the following countable equivalence relations are ultrahomogeneous: $\bigcup _{\o }F_n$ (indivisible iff $n=1$); $\bigcup _{n}F_\o$ (indivisible iff $n=1$) and
$\bigcup _{\o }F_\o $ (the $\o$-homogeneous-universal equivalence relation, indivisible of course).
\\
The same picture is obtained for

- Disconnected countable ultrahomogeneous graphs,
which are (by the well known classification of Lachlan and Woodrow) of the form
$\bigcup _m K_n$, where $mn=\o$ (the disjoint union of $m$-many complete graphs of size $n$);

- Countable posets of the form $\bigcup _m L_n$, where $mn=\o$ (the disjoint union of $m$-many copies of the ordinal $n\in [1,\o]$).

We note that the relational structures observed in this section are disconnected but taking their complements we obtain
connected structures with the same posets $\la \P (\X ), \subset\ra$ and $\sq\la \P (\X ), \subset\ra$. For example, the complement of $\bigcup _m F_n$ is  the graph-theoretic complement of the graph $\bigcup _m K_n$.

\end{document}